\documentclass[11pt]{amsart}
\usepackage{tikz,amsthm,amsmath,amstext,amssymb,amscd,epsfig,euscript, mathrsfs, dsfont,pspicture,multicol,graphpap,graphics,graphicx,times,enumerate,subfig,sidecap,wrapfig,color}

\usepackage{amsmath}
\usepackage{amssymb}

\newcommand{\R}{{\mathbb R}}       
\newcommand{\N}{{\mathbb N}}       %
\newcommand{\Z}{{\mathbb Z}}       

\newcommand{\diam}{\mathop{\rm diam}}

\newcommand{\vv}{{\vspace{2mm}}}

\newcommand{\noi}{\noindent}

\def\XXint#1#2#3{{\setbox0=\hbox{$#1{#2#3}{\int}$ }
\vcenter{\hbox{$#2#3$ }}\kern-.58\wd0}}

\newtheorem{theorem}{Theorem}[section]

\newtheorem*{lemma*}{Lemma}
\newtheorem*{theorem*}{Theorem}

\theoremstyle{definition}

\theoremstyle{remark}
\newtheorem{rem}[theorem]{\bf Remark}

\numberwithin{equation}{section}

\newcommand{\RRem}{\begin{rem}}
\newcommand{\erem}{\end{rem}}

\makeatletter
\def\@tocline#1#2#3#4#5#6#7{\relax
  \ifnum #1>\c@tocdepth 
  \else
    \par \addpenalty\@secpenalty\addvspace{#2}%
    \begingroup \hyphenpenalty\@M
    \@ifempty{#4}{%
      \@tempdima\csname r@tocindent\number#1\endcsname\relax
    }{%
      \@tempdima#4\relax
    }%
    \parindent\z@ \leftskip#3\relax \advance\leftskip\@tempdima\relax
    \rightskip\@pnumwidth plus4em \parfillskip-\@pnumwidth
    #5\leavevmode\hskip-\@tempdima
      \ifcase #1
       \or\or \hskip 1em \or \hskip 2em \else \hskip 3em \fi%
      #6\nobreak\relax
    \dotfill\hbox to\@pnumwidth{\@tocpagenum{#7}}\par
    \nobreak
    \endgroup
  \fi}
\makeatother

\def\cB{{\mathscr{B}}}

\def\cD{{\mathscr{D}}}
\def\cE{{\mathscr{E}}}
\def\cF{{\mathscr{F}}}

\def\cH{{\mathscr{H}}}

\def\cR{{\mathscr{R}}}
\def\cS{{\mathscr{S}}}

\begin{document}

\title[Extrapolation of Carleson measures] {On Extrapolation of Carleson measures}

\author[Garnett]{John Garnett}

\address{John Garnett\\
Department of Mathematics, 6363 Mathematical Sciences Building \\University of California
at Los Angeles, 520 Portola Plaza, Los Angeles, California 90095-1555.
}
\email{jbg@math.ucla.edu}


\begin{abstract}

Extrapolation of Carleson measures is a technique originated by [LM] which has become an essential step for several important results.  See for example [HL], [AHLT], [AHMTT] and more recently [HM], [HMM], and [HMMTZ].  Here we recast the proof of  the extrapolation theorems in [AHLT] and [HMM]  using the stopping 
time - generations language from [C1], [C2] or [CG]  but no new idea. We give three equivalent versions of the theorem having the same proof.   The first is for dyadic cubes in $[0,1]^d,$ the second is an abstract formulation of the result on Ahlfors regular sets from [HM], and the third is the version in [AHLT] which follows from the other two.

\end{abstract}


\maketitle



\section{The dyadic case}

In $\R^d, ~ d \geq 1$, let $\cD$ denote the set of closed dyadic cubes 

$$Q = \bigcup_{k=1}^d \{j_k2^{-n} \leq x_k \leq  (j_k+1)2^{-n}\} \subset Q_0 = [0,1]^d, n \in \N, j_k \in \Z.$$ 

\noi For $Q \in \cD$ write $\ell(Q) = 2^{-n}$ for its side length and $|Q| = 2^{-nd}$ for its measure and     
define 

$$Q^* = Q \times (0,\ell(Q)] \subset \R^{d+1}_+,$$
 
\vv 

$$T(Q) = Q \times (\ell(Q)/2, \ell(Q)] = 
{Q}^* \setminus \cup\{Q_1^*: Q_1 \subsetneqq Q\},$$ 

\noi and 

$$\cD(Q) = \{Q' \in \cD: Q' \subset Q\}.$$

\vv

A Borel measure $\mu$ on $Q_0^*$ is a {\bf Carleson measure} if 

\begin{equation}\label{(1.1)}
C_1(\mu) =\sup_{Q \in \cD} \frac{\mu(Q^*)}{|Q|} < \infty.
\end{equation}

\noi Let $\nu$ be another Borel measure on $Q_0^*$ satisfying 

\begin{equation}\label{(1.2)}
C_2(\nu) = \sup_{Q \in \cD} \frac{\nu(T(Q))}{|Q|} < \infty,
\end{equation}  

\noi a condition  clearly necessary for $\nu$ to be a Carleson measure.  

\vv

\noi {\bf  {Theorem 1.1:}}  Let $\mu$ and $\nu$ be Borel measures on $Q_0^*$ satisfying  (1.1) and (1.2) respectively.  Assume there exist constants $\delta > 0$ and $C$ such that 

\begin{equation}\label{(1.3)}
\nu(Q^* \setminus \bigcup_{\cF} Q'^*)  \leq C |Q|
\end{equation}

\noi whenever $Q \in \cD$ and $\cF \subset \cD(Q)$  is a set of subcubes of $Q$ with  disjoint interiors for which 

\begin{equation}\label{(1.4)}
\sup_{Q' \in \cD(Q)} 
\frac {\mu(Q'^* \setminus \bigcup_{\cF} Q''^*)}{|Q'|} \leq \delta.
\end{equation}

\noi Then $\nu$ is a Carleson measure with constant 

\begin{equation}\label{(1.5)}
C_1(\nu) \leq \biggl(C_2(\nu) + \frac{2^dC}{2^d-1}\biggr)\frac{C_1(\mu)}{\delta}.
\end{equation}

\vv

Careful reading will show that the proof below is little more than a reformulation of the arguments in Sections 7 and 9 of [HM].

\vv

\noi {\bf Proof:} When proving Theorem 1.1. we can assume 
$E = \bigcup_{\cD}  \partial Q^*$ satisfies $\mu(E) = \nu(E) =0$ by pushing some mass off each $\partial Q^*$ and thus we can treat all rectangles as open or closed sets when taking unions or intersections.

\vv

Fix $\delta > 0$ for which (1.4) implies (1.3).  For each $Q \in \cD$ we will define a subset $U(Q) \subset Q^*$  and a pairwise disjoint family $G_1(Q) \subset \cD(Q)$ so that 

\begin{equation} \label{(1.6)}
U(Q) = Q^* \setminus \bigcup_{G_1(Q)} Q'^*.
\end{equation}

\noi Further define by induction 

$$G_n(Q) = \bigcup \{G_1(Q'): Q' \in G_{n -1}(Q)\}$$ 

\noi and $G_0(Q) = \{Q\}.$   Then the family 
$\bigcup_{n=0}^{\infty} \{U(Q'): Q' \in G_n(Q)\}$ 
is pairwise disjoint and 

\begin{equation}\label{(1.7)}
\nu(Q^*) = \sum_{n=0}^{\infty}\sum_{G_n(Q)} \nu(U(Q')).
\end{equation}

\vv

\noi For each $Q$ and $n \geq 1$ write

$$\cE_n(Q) = \{Q_n \subset Q: \ell(Q_n) = 2^{-n} \ell(Q)\}.$$

\noi Also define

\begin{equation}\label{(1.8)}
\cB = \{Q:  \mu(T(Q)) \geq \delta |Q|\}.
\end{equation}

\vv

 Fix $Q \in \cD$.  For $Q' \subset Q$ the definitions of $U(Q')$ and $G_1(Q')$ depend on whether or not $Q' \in \cB.$ 
When $Q \in \cB$ take $G_1(Q) = \cE_1(Q)$ and $U(Q) = T(Q).$  Then for any $Q$  (1.2), (1.8) and disjointness yield 

\begin{equation}\label{(1.9)}
\sum_{n=0}^{\infty} \sum_{\cB \cap G_n(Q)} \nu(U(Q')) 
\leq \frac {C_2}{\delta}
\sum_{\cB \cap \cD(Q)} \mu(T(Q'))  \leq  \frac {C_1 C_2}{\delta} |Q|, 
\end{equation}

\noi so that when estimating (1.7) we need only  study $Q' \notin \cB.$

\vv

When   $Q \in \cD \setminus \cB$ we define  by induction sets $\cF_n(Q) \subset \bigcup_{k=1}^n \cE_k(Q)$ so that 

$$\cF_n(Q) \subset \cF_{n+1}(Q)$$

\noi  and so that (1.4) holds  for $Q$ and 

$$\widetilde{\cF_n}(Q) = 
\cF_n(Q) \cup \cR_{n+1}(Q)$$

\noi where 

$$\cR_{n+1}(Q) = \cE_{n+1}(Q)  \setminus \bigcup_{Q' \in \cF_n(Q)} 
\cD(Q').$$

\noi Then (1.4) will hold for $Q$ and

$$\cF(Q) = \bigcup_n \cF_n(Q).$$

\vv  

\noi {\bf {Step I:}}  $n =1.$ Include $\cB \cap \cE_1(Q)  \subset \cF_1(Q)$ and then consider the family of subsets
$\cS \subset \cE_1(Q) \setminus \cB$ such that (1.4) holds  for $Q$ 
and $\widetilde \cF_{\cS}$ where 
$\cF_{\cS} = \cS \cup (\cE_1(Q) \setminus \cB)$. This  family is non-empty because it contains 
$\cE_1(Q) \setminus \cB$ since $Q \notin \cB.$
Order this set family by inclusion,  let $\cS_1(Q)$ be a minimal element and define  

$$\cF_1(Q) = \cS_1(Q) \cup (\cB \cap \cE_1(Q)).$$ 

\noi  Then (1.4) holds for  
$\widetilde \cF_1(Q).$  

\vv  

\noi {\bf {Step II:}}   Now assume 
$n \geq 1$ and $\cF_n(Q)$ has been constructed. If $\bigcup_{\cF_n(Q)} Q' = Q$ stop the construction for $Q$.   If not, include 
$\cB \cap \cR_n(Q)  \subset \cF_{n+1}(Q)$
\noi and consider the family of subsets 
$\cS \subset  \cR_n(Q) \setminus \cB $ such that  when

\begin{equation}\label{(1.10)}
\cF_{\cS}  = \cF_n(Q) \cup \cS \cup (\cB \cap \cR_{n+1}(Q))
\end{equation}

\noi (1.4) holds for $Q$ and 

$$\widetilde \cF_{\cS} = \cF_{\cS} \cup \bigl(\cE_{n+2}(Q) \setminus \bigcup_{Q' \in \cF_{\cS}} \cD(Q')\bigr).$$

\noi This family is non-empty because by induction (1.4) holds  for $Q$ and
 $\widetilde \cF_{\cS}$ when $\cS = \cR_{n+1}(Q).$ Order this finite set family by inclusion, let
$\cS_{n+1}(Q)$ be a minimal member, and define
$\cF_{n+1}(Q)$ by (1.6) with $\cS = \cS_{n+1}(Q).$   Then (1.4) holds for $Q$ and $\widetilde{\cF_{n+1}}(Q).$ 
Moreover we have:
\vv 

\noi {\bf {Lemma 1.2:}} If $Q' \in \cF(Q) \setminus \cB$ there exists 
$\widetilde Q'$ such that $Q' \subset \widetilde Q' \subset Q$ and  

\begin{equation}\label{(1.11)}
\mu(\widetilde Q'^* \setminus \bigcup_{\cF(Q)} Q''^*) \geq \biggl(1 -\frac{1}{2^d}\biggr)
\delta|\widetilde Q'|.
\end{equation}

\vv 

\noi {\bf Proof:}  We have $Q' \in \cF_n$ for some first $n$ and by the minimality 
of $\cS_n$  there exists  $\widetilde Q' \supsetneqq Q'$ such that 

$$\mu\bigl(\widetilde Q'^* \setminus \bigcup \{Q''^*: Q'' \in \cF_n(Q) \cup\cE_{n+1}(Q)\}\bigr) + \mu(T(Q')) > \delta |\widetilde Q'|$$

\noi while $\mu(T(Q')) < \frac{\delta}{2^d} |\widetilde Q'|$ since $Q' \notin \cB.$  Therefore (1.11) holds. 

\vv

When $Q' \notin \cB$, (1.4) holds for $Q'$ and $\cF(Q')$ so that by (1.3)
$\nu(U(Q')) \leq C |Q'|$ and  therefore 

\begin{equation}\label{(1.12}
\sum_n \sum_{G_n(Q) \setminus \cB} \nu(U(Q')) 
\leq C \sum_{n=0}^{\infty} \sum_{G_n(Q) \setminus \cB} |Q'|.
\end{equation}

\noi By  Lemma 1.2 we can cover  $\bigcup_{\cF(Q) \setminus \cB} Q'$ by the pairwise disjoint family $\cH(Q)$ of maximal dyadic cubes $\widetilde Q'$ that satisfy (1.11).  Then since both  families $\cH(Q)$ and $\cF(Q) \setminus \cB$ are pairwise disjoint. (1.11) yields

\begin{equation}\label{(1.13)}
\sum_{\cF(Q) \setminus \cB} |Q'|  
 \leq \sum_{\cH(Q)} |\widetilde Q'| \leq \frac {2^d}{2^d-1} \frac {1}{\delta} \mu(U(Q)),
\end{equation}

\noi and (1.13),  (1.12) and (1.1) give

\begin{equation}\label{(1.14)}
\sum_{n =1}^{\infty} \sum_{G_n(Q) \setminus \cB} \nu(Q') 
\leq  \frac{2^d}{2^d-1}\frac{CC_1}{\delta}|Q|,
\end{equation}

\vv
\noi and together (1.14)  and (1.9) establish the estimate (1.5) and prove Theorem 1.1. 

\vv
\vv

\section {a general case}

Let $d$ be a positive integer and let $(X,\rho)$ be a metric space on which there exists a positive Borel measure $\sigma $ that is $d$-Ahlfors regular: there exists $c_1 > 0$ such that 

\begin{equation}\label{(2.1)}
\frac{1}{c_1} R^d \leq \sigma(B(x,R) \leq c_1 R^d
\end{equation} 

\noi for all $x \in X$ and all $0 < R \leq \diam(X).$  We assume for 
simplicity that $X$ is compact and  $\sigma(X) =1.$  Following Christ [Ch], there exists a positive integer $N$ and a  family

$$\cD = \bigcup_{j = 0}^{\infty} \cD_j$$

\noi of Borel subsets of $X$ satisfying  (2.2) - (2.6) below:


\begin{equation}\label{(2.2)}
\diam{Q} \sim 2^{-Nj}  ~~~{\rm {if}}~~~ Q \in \cD_j;
\end{equation}

\vv

\begin{equation}\label{(2.3)}
X = \bigcup_{{\cD}_j}Q, ~~~{\rm {for ~all}}~~~ j;
\end{equation}

\vv

\begin{equation}\label{(2.4)}
Q \cap  Q'  = \emptyset ~~~ {\rm {if}} ~~~ Q,~ Q'\in {\cD}_j ~~~{\rm {and}}~~~ Q' \neq Q;
\end{equation}

\vv

\begin{equation}\label{(2.5)}
{\rm {if}}~~ j < k, ~~~ Q_j \in \cD_j ~~{\rm {and}} ~Q_k \in \cD_k, ~~~ {\rm {then}} ~~~  Q_k \subset Q_j ~~~{\rm {or}}~~~ 
Q_k \cap Q_j = \emptyset. 
\end{equation}

\vv

\noi There exists constant $c_0$ such that for all $Q \in \cD$ there exists $x_Q \in Q$ with  

\begin{equation}\label{(2.6)}
B(x_Q,c_0\ell(Q)) \cap \partial \Omega \subset Q.
\end{equation}

Note that by (2.1), (2.2) and (2.6) there is a constant $c_2$ so that for all $Q \in \cD$,

\begin{equation}\label{(2.7)}
\frac{1}{c_2} \ell(Q)^d \leq \sigma (Q) \leq c_2 \ell(Q)^d.
\end{equation}

Now let $\mu$ and $\nu$ be positive discrete measures on the countable 
set $\cD$, so that there exist $\alpha_Q \geq 0$ and $\beta_Q \geq 0$ such that for any $\cE \subset \cD$

$$\mu(\cE) = \sum_{Q \in \cE} \alpha_Q$$

\noi and 

$$\nu(\cE) = \sum_{Q \in \cE} \beta_Q.$$

\noi Analogous to (1.1) and (1.2) we assume $\mu$ is a discrete Carleson measure,

\begin{equation}\label{(2.8)}
C_1(\mu) = \sup_{Q \in \cD}\frac{\mu(Q^*)}{\sigma(Q)} < \infty,
\end{equation} 

\noi where we write $Q^* = \{Q': Q' \subseteq Q\}$, and we assume 

\begin{equation}\label{(2.9)}
C_2(\nu)  \sup_{Q \in \cD} \frac{\beta_Q}{\sigma(Q)} < \infty,
\end{equation}

\noi which is necessary for $\nu$ to be a discrete Carleson measure. Then we have the following abstract version of Theorem 1.1:

\vv

\noi {\bf  {Theorem 2.1:}}  Let $\mu$ and $\nu$ be measures on $\cD$ satisfying  (2.8) and (2.9) respectively.  Assume there exist constants $\delta > 0$ and $C$ such that 

\begin{equation}\label{(2.10)}
\nu(Q^* \setminus \bigcup_{\cF} Q'^*)  \leq C \sigma(Q)
\end{equation}

\noi whenever $Q \in \cD$ and $\cF \subset \cD(Q)$  is a set of subcubes of $Q$ for which 

\begin{equation}\label{(2.11)}
\sup_{Q' \subset Q} 
\frac {\mu(Q'^* \setminus \bigcup_{\cF} Q''^*)}{\sigma(Q')} \leq \delta. \end{equation}

\noi Then $\nu$ is a Carleson measure with constant 
$C_1(\nu)$ depending only on $c_2, C, \delta, C_1(\mu)$ and  $C_2(\nu).$

\vv  

But for the choices of constants the proof of Theorem 2.1 is a repetition of the proof of Theorem 1.1 and we omit the details. 

 Theorem 2.1 can be stated more abstractly as a result about the tree $(\cD, \subset)$ with transition  probabilities $\frac{\sigma(Q')}{\sigma (Q)}$ for  $Q' \subset Q.$ 
However for the important applications in [HM], [HMM] and [HMMTZ]
$\Omega$ is a domain in $\R^{d+1}$,  $E = \partial \Omega$ is uniformly rectifiable,
and $\alpha_Q$ and $\beta_Q$ depend critically on the uniform rectifiabity properties of $E$ and on the elliptic differential equation whose solutions are being estimated.  Thus the real difficulty lies in finding suitable functions
$\alpha_Q$ and $\beta_Q.$


\vv 
\vv

\section{a third version}

The ``extrapolation lemma" in [AHLT]  assumes $\mu$ is a Carleson measure in 
$R^{d+1}_+$, i. e. $\mu$ satisfies (1.1) for all cubes $Q \subset \R^d$ and  $\nu$ is another Borel measure on $R^{d+1}_{+}$ satisfying 

\begin{equation}\label{(3.1)}
\nu(T(Q)) \leq C_2 |Q|,   
\end{equation}

\noi where now $T(Q) = Q \times [\frac{\ell(Q)}{2},\ell(Q)] \subset \R^{d+1}_{+}.$

But instead of a cube family $\cF$ one now works with a nonnegative Lipschitz function $\psi$ such that 

\begin{equation}\label{(3.2)}
||\nabla \psi||_{\infty} \leq 1
\end{equation}  

\noi and the region

$$\Omega_{\psi} = \{(x,t) \in R^d \times (0,\infty): t \geq \psi(x)\}$$

\vv
\noi and instead of (1.4) one tests $\mu(T_Q \setminus \Omega_{\psi})$ where $T_Q$ is the tent

$$T_Q = \{(x,t): x \in Q, 0 \leq t \leq {\rm {dist}}(x, \R^d \setminus Q)\}.$$

\noi {\bf {Theorem 3.1:}} Let $\mu$ and $\nu$ be Borel measures on $R^{d+1}_{+}$ satisfying (1.1) and (3.1) respectively. Assume there is are constants $\delta > 0$ and $C > 0$ such that whenever $Q \subset \R^d$ is a cube and $\psi:Q \to \R$ is a nonnegative Lipschitz function satisfying (3.2)

\begin{equation}\label{(3.3)}
\nu(Q^* \cap \Omega_{\psi}) \leq C |Q|
\end{equation}

holds provided

\begin{equation}\label{(3.4)}
\sup_{Q' \subset Q} \frac{\mu(T_{Q'} \cap \Omega_{\psi})}{|Q'|} \leq \delta
\end{equation}

\noi where the supremum is taken over the decomposition of $Q$ into cubes having side 
$2^{-n} \ell(Q),  n \in \N.$ . Then $\nu$ is a Carleson measure,  i.e. for all $Q$ 

\begin{equation}\label{(3.5)}
\nu(Q^*) \leq C |Q|
\end{equation}

\vv

\noi {\bf {Proof:}} First note that (3.1) implies for all $Q$ 

\begin{equation}\label{(3.6)}
\nu(Q^* \setminus T_Q) \leq (C_2 +d)|Q|
\end{equation}

\noindent so that to prove (3.4) for a fixed cube $Q$ we may assume 

\begin{equation}\label{(3.7)}
\nu(Q* \setminus T_Q) =0
\end{equation}

\noi and as before we can assume $\nu(\partial T_Q' \cup T(Q')) = 0$ for all $Q' \subset Q$. For each $Q'$ write $p_{Q'} = (c(Q'), \frac{\ell(Q'}{2})$ for the center of $Q^*$ which is also the vertex of $T_{Q'}$  and define the discrete measures 

$$\tilde \nu = \sum_{Q'\subsetneqq Q} \nu(T(Q'))\delta_{p_{Q'}}, ~~{\rm {and}} ~~\tilde \mu = 
\sum_{Q'\subsetneqq Q} \mu(T(Q'))\delta_{p_{Q'}},$$

\noi where $\delta_p$ is the unit point mass at $p$.   Then by (3.7)  

\begin{equation}\label{(3.8)}
\tilde \nu(Q'^*) = \nu(Q'^*) ~~{\rm {and}}~~ \tilde \mu(Q'^*) = \mu(Q'^*)
\end{equation} 

\noi for all $Q' \subset Q$.  

\vv
To prove Theorem 3.1 we show that if (3.4) $\Longrightarrow$ (3.3) for $\mu$ and $\nu$ then
(1.4) $\Longrightarrow$ (1.3) for $\tilde \mu$ and $\tilde \nu$, so that by (3.8) and Theorem 1.1, (3.5) holds  
for $\nu.$ To that end let $\cF = \{Q_1, ....\}$ be a pairwise disjoint family of subcubes of $Q$ satisfying (1.4) with $\tilde \nu.$ Set

$$\psi_{\cF}(x) = \sum_{\cF} \frac{2}{\ell(Q_j)}\chi_{Q_j}(x) {\rm {dist}}(x,\partial Q_j).$$  

\noi Then 

$$Q^* \cap \Omega_{\psi} = \bigl(Q^* \setminus \bigcup_{\cF} Q_j^*\bigr) \cup \bigcup_{\cF} 
(Q_j^* \setminus T_{Q_j})$$

\noi so that by (1.4) and (3.8) we obtain (3.4) for $\mu$ and (3.3) for $\nu$ and therefore (1.3) for 
$\tilde \nu$, again by (3.8).

\vv

\end{document}